\RequirePackage{fix-cm}
\documentclass[smallextended]{svjour3}       
\smartqed  
\usepackage{graphicx}
\usepackage{amsmath, amssymb}
\usepackage{multirow}
\begin{document}

\title{Probability Distribution-free General Scenario Programming
}
\subtitle{Applications in Smart Grid Optimization under Both Exogenous and Endogenous Uncertainty}


\author{Qifeng Li          
}


\institute{Qifeng Li \at
              Department of Electrical and Computer Engineering, University of Central Florida, Orlando, FL 32816, USA\\
              \email{@ucf.edu}           
}

\date{Received: date / Accepted: date}

\maketitle

\begin{abstract}
This paper presents a novel solution paradigm of general optimization under both exogenous and endogenous uncertainties. This solution paradigm consists of a probability distribution (PD)-free method of obtaining deterministic equivalents and an innovative approach of scenario reduction. First, dislike the existing methods that use scenarios sampled from pre-known PD functions, the PD-free method uses historical measurements of uncertain variables as input to convert the logical models into a type of deterministic equivalents called General Scenario Program (GSP). Our contributions to the PD-free deterministic equivalent construction reside in generalization (making it applicable to general optimization under uncertainty rather than just chance-constrained optimization) and extension (enabling it to the problems under endogenous uncertainty via developing an iterative and a non-iterative frameworks). Second, this paper reveals some unknown properties of the PD-free deterministic equivalent construction, such as the characteristics of active scenarios and repeated scenarios. Base on this discoveries, we propose a concept and methods of strategic scenario selection which can effectively reduce the required number of scenarios as demonstrated in both mathematical analysis and numerical experiments. Numerical experiments are conducted on two typical smart grid optimization problems under exogenous and endogenous uncertainties.

\end{abstract}

\section{Introduction}
\label{intro}
This paper considers the following logical model of optimization under uncertainty (OU):
\begin{subequations} \label{OU}
\begin{align}
\text{(OU)}\quad\quad \min_{x} \quad & f_0(x) + \mathbb{E}_{w}[f_1(x,w)] \label{obj_OU} \\
 \mathrm{s.t.} \quad  & \mathbb{P}_{w}[g(x,w) \le 0] \ge 1-\alpha, \label{constr_OU} 
\end{align}
\end{subequations}
where $x \in \mathbb{R}^n$ is the vector of (continuous, integer, or mixed continuous and integer) decision variables, $f_0:\mathbb{R}^n \rightarrow \mathbb{R}$, $f_1:\mathbb{R}^{n+r} \rightarrow \mathbb{R}$, and $g:\mathbb{R}^{n+r} \rightarrow \mathbb{R}^m$ are general functions, $\mathbb{P}$ and $\mathbb{E}$ represent probability and expectation respectively, and $\alpha \in \mathbb{R}^m$ represent the allowed probabilities of constraint violation. For simplicity but without loss of generality, only inequality constraints are considered in (\ref{constr_OU}) since there are explicit and inexplicit methods of equivalently representing equations as inequalities. The vector of uncertainties $w =[u,v]^{\rm T} \in \mathbb{R}^r$ consists of exogenous uncertainties (ExU) $u\in \mathbb{R}^{r_1}$ and endogenous uncertainties (EnU) $v\in \mathbb{R}^{r_2}$, which follow certain probability distributions (PD):
\begin{subequations}
\begin{align}
    & \mathbb{P}[u_i]=p_i(u_i), \;u_i \in \mathcal{U},\; i=1,\ldots,r_1 \label{PDx}\\
    & \mathbb{P}[v_j]=q_j(v_j,x), \;v_j \in \mathcal{V},\; j=1,\ldots,r_2 \label{PDn}
\end{align}
\end{subequations}
where ($r_1+r_2=r$), $\mathcal{U}$ and $\mathcal{V}$ are the uncertain sets. Here, ExU and EnU mean the decision-independent and decision-dependent uncertainties respectively \cite{nohadani2018optimization}.

Formulation (OU) is a general logical model of optimization under uncertainty. Given that (\ref{constr_OU}) is a set of chance constraints, (OU) represents the chance-constrained optimization (CCO \cite{li2008chance}) if $f_1=0$. When $\alpha=0$, (OU) reduces to a stochastic optimization (SO \cite{birge2011introduction}) and further reduces to a robust optimization (RO \cite{bertsimas2011theory}) if $f_1=0$. However, problem (OU) can not be directly solved by mature optimization algorithms (e.g., interior-point, quasi-Newton, and conjugate gradient \cite{nocedal2006numerical}) or  computer solvers that are based on these algorithms (e.g., Knitro, Gurobi, and Mosek \cite{anand2017comparative}). In the solution process, an inevitable step is converting the logical model (OU) into its deterministic equivalents, i.e. linear, nonlinear, or integer program. Inspired by the widely adopted term ``convexification" in the optimization field, which refers to the process of converting or approximating non-convex problems into convex ones, we define a verb ``determinisfy" to refer to convert a logical model into its deterministic equivalent, and the term ``deterministication" to refer to the process of determinisfying.

The existing deterministication methods assume that the PD functions are known a prior and determinisfy the logical model (OU) using scenarios of the uncertain variables sampled from the PD functions \cite{birge2011introduction}. However, such PD-based deterministication methods creates an ``exponential dilemma", i.e. the problem sizes of the resulting deterministic equivalents increase \textit{exponentially} as the number of uncertain variables grows. For instance, $a$ samples for each uncertain variable out of $b$ will result in $a^b$ scenarios. The authors of \cite{campi2008exact,calafiore2010random} used randomly selected scenarios to convert a convex CCO under only ExU into a uncertain/random convex program which is a deterministic equivalent candidate of the original CCO. They also provided a closed-form relation between the numbers of input scenarios and the degree of ``equivalent." In the past decade, researchers has attempted to extend this PD-free method to continuous nonconvex \cite{campi2018general} and mixed-integer convex cases \cite{calafiore2012mixed,esfahani2014performance}. However, we found that the numbers of needed scenarios determined in these literature are unnecessarily high, so that the applicability of this approach to large engineering systems is very limited. Moreover, EnU is not considered in these existing research.

Based on the above-mentioned PD-free deterministication method, this paper develops a novel technical path for effectively solving a more general problem, i.e. (OU). On this research path, we aim at achieving a high equivalence using as small numbers of scenarios as possible rather than finding the exact closed-form formulation of the relation between the numbers of input scenarios and the degree of ``equivalent.". As a result, this deterministication method can be applied to general optimization under uncertainty. More importantly, we will enable this method for determinisfying problems with EnU, which has not been investigated in existing research.  

\section{PD-Free Deterministication under Exogenous Uncertainty}
\label{sec:1}

For simplicity, we start our discussions with the logical model (OU) under only ExU, where $w$ reduces to only $u$. This section aims at establishing a methodoglogy of PD-Free deterministication by generalizing the existing findings on scenario optimizations.

\subsection{General Scenario Program}
Denoting the (OU) under ExU as (OU-ExU), the developed PD-Free Deterministication method determinisfies the logical model (OU-ExU) into the following deterministic equivalent: 
\begin{subequations} \label{GSP}
\begin{align}
\text{(GSP-ExU)}\quad\quad \min_{x} \quad & f_0(x) +\frac{1}{N}\sum_i^N f_1(x,u^{(i)})   \label{obj_GSP} \\
 \mathrm{s.t.} \quad  &g(x,u^{(i)}) \le 0,\;   (u^{(i)}  \in \mathcal{U},\;i=1,\ldots,N)\label{constr_GSP} 
\end{align} 
\end{subequations}
which is called general scenario program (GSP) in this paper. In the (GSP-ExU) formulation, $\mathcal{U}$ is a big enough set of historical measurements of $u$. Since $u^{(i)}$ ($i=1,\ldots,N$) are $N$ scenarios \textit{randomly} picked up from $\mathcal{U}$, the transformation from (OU-ExU) to (GSP-ExU) is free from PD functions. Let $x^*_N$ be the optimal solution of (GSP-ExU) with $N$ input scenarios, we have the following definition.
\smallskip

\noindent
{\bf Definition 1}: A confidence factor $\epsilon$ ($0 \le \epsilon \le 1$) is defined to capture the probability that $x^*_N$ is optimal to (OU-ExU), namely $\epsilon=\min \{\epsilon_1,\epsilon_2\}$ and
\begin{subequations}
\begin{align}
  \epsilon_1 =  &\mathbb{P}_N[\mathbb{P}_{u}[g(x^*_N,u) \le 0] \ge 1-\alpha] \\
 \epsilon_2 = &\mathbb{P}_N[f_0(x^*_N)+\mathbb{E}_{u}[f_1(x^*_N,u) \le f_0(\tilde{x})+\mathbb{E}_{u}[f_1(\tilde{x},u)]\\
 &\text{if}\; \mathbb{P}_{u}[g(\tilde{x},u) \le 0] \ge 1-\alpha, \nonumber
\end{align}
\end{subequations}
where $u \in \mathcal{U}$. 

The confidence factor $\epsilon$ can be interpreted as the degree that the deterministic (GSP-ExU) is equivalent to (OU-ExU) or the accuracy of using the deterministic (GSP-ExU) to approximate (OU-ExU). Since $\epsilon$ represents a desired confidence level, it is a given value for example 0.99.
\smallskip

\subsection{Existing Findings on Chance-constrained Optimization under ExU}

This subsection considers a subset of (OU-ExU):
\begin{subequations} \label{CCOx}
\begin{align}
\text{(CCO-ExU)}\quad\quad \min_{x} \quad & f_0(x) \label{obj_CRP} \\
 \mathrm{s.t.} \quad  & \mathbb{P}_{u}[g(x,u) \le 0] \ge 1-\alpha, \label{constr_CRP} 
\end{align}
\end{subequations}
which is a chance-constrained optimization (CCO) under ExU. Following (GSP-ExU), the corresponding deterministic equivalent of (CCO-ExU) in the form of GSP is given as:
\begin{subequations} \label{CCO-GSP}
\begin{align}
\text{(GSP-CCO-ExU)}\quad\quad \min_{x} \quad & f_0(x)   \label{obj_CSP} \\
 \mathrm{s.t.} \quad  &g(x,u^{(i)}) \le 0.\;   (u^{(i)}  \in \mathcal{U},\;i=1,\ldots,N)\label{constr_CSP} 
\end{align} 
\end{subequations}
The existing research mainly covers three special cases of (GSP-CCO-ExU) as tabulated in Table 1. The findings of the existing research are summarized in the following proposition.
\smallskip

\begin{table}[h]
\centering
\caption{The three special cases of (CCO-ExU) in existing research}
\begin{tabular}{|c|c|c|c|}
\hline
                            & \textbf{Case 1}  & \textbf{Case 2}  & \textbf{Case 3}   \\ \hline
\textbf{Objective function} & $f_0(x)=c^{\rm T}x$               & $f_0(x)=c^{\rm T}x$                & $f_0$ is any function \\ \hline
\textbf{Constraints}        & $g$ is convex on $x$ & $g$ is convex on $x$ & $g$ is any function \\ \hline
\textbf{Variables}          & continous        & mixed-integer    & continous         \\ \hline
\end{tabular}
\end{table}

\noindent
{\bf Proposition 1} on (GSP-CCO-ExU): $x^*_N$ \textit{is optimal to} (CCO-ExU) \textit{with a confidence level of $\epsilon$ if,}

\noindent
\text{1.}\; \textit{for Case 1:}
\begin{equation}
 N \ge \frac{e}{e-1}\frac{1}{\alpha}(n+\ln{\frac{1}{1-\epsilon}})\; \label{N1}
\end{equation}

\noindent
\text{2.}\;\textit{for Case 2:}
\begin{equation}
N  \ge \underset{N^{\min}}{\arg}\left\{\sum_{k=0}^{(n_1+1)2^{n_2}} [\frac{N^{\min}!}{k!(N^{\min}-k)!}\alpha^k(1-\alpha)^{N^{\min}-k}] =1-\epsilon \right\}\;  \label{N2}
\end{equation}

\noindent
\text{3.}\;\textit{for Case 3:}\; 
\begin{align}
    N  \ge &\underset{N^{\min}}{\arg} \{\sum_{k=0}^{N^{\min}} [\frac{N^{\min}!}{k!(N^{\min}-k)!}(1-\frac{1}{N^{\min}-k}\log \frac{1}{\beta} \nonumber \\
&-\frac{1}{N^{\min}-k}\log \frac{N^{\min}!N^{\min}}{k!(N^{\min}-k)!})] =1-\epsilon \}, \label{N3}
\end{align}
\textit{where $e$ is Euler's number, $n_1$ and $n_2$ are the numbers of continuous and integer variables, and $N^{\min}$ is the solution of the equation inside the curly brackets.}
\smallskip

\textit{Proof:} Find the proof of Case 1 in \cite{campi2008exact,calafiore2010random}, the proof of Case 2 in \cite{calafiore2012mixed}, and the proof of Case 3 in \cite{campi2018general}.

\hfill$\square$


\subsection{A Hypothesis on the PD-free Deterministication of (OU-ExU)}
This subsection introduces our hypothesis on the PD-free Deterministication we established in Subsection 2.1.
\smallskip

\noindent
{\bf Hypothesis}: \textit{Let $x^*_N$ be the optimal solution of} (GSP-ExU) \textit{with $N$ input scenarios, then:}
\vspace{-3pt}
\begin{enumerate}
    \item \textit{There exists a number $N^{\min}$ that, if $N \ge N^{\min}$, $x^*_N$ is also optimal to} (OU-ExU) \textit{with a confidence factor }$\epsilon$.
\item \textit{For a specific case, the relation among} $N^{\min}$, $\epsilon$, $\alpha$, $n$ \textit{and} $r$ \textit{can be formulated by a closed-form expression} $N^{\min}=h(\epsilon, \alpha, n, r)$ \textit{where:}
\begin{equation}
    \frac{\partial h}{\partial \epsilon},\,\frac{\partial h}{\partial n},\, \frac{\partial h}{\partial r} > 0\;\text{and}\;\frac{\partial h}{\partial \alpha}<0. \label{hypo}
\end{equation}
\end{enumerate}
\smallskip

\textit{A qualitative explanation:} If a scenario $u^{(j)}$'s probability is $\mathbb{P}[u^{(j)}]$, it will \textit{most likely} appear ($N*\mathbb{P}[u^{(j)}]$) times in the selected scenario set $\mathcal{N}=\{u^{(i)}, i=1,\ldots,N\}$ according to the basic theory in statistics. If ($N*\mathbb{P}[u^{(j)}] \ge 1$), $u^{(j)}$ is most likely selected and inputted to (\ref{GSP}), such that the optimal solution $x^*_N$ is feasible to $u^{(j)}$. Otherwise, $u^{(j)}$ is most likely not selected and $x^*_N$ is not necessarily feasible to $u^{(j)}$. In other words, if $x^*_N$ is required to be feasible to all scenarios whose probability is equal or bigger than $\alpha$, one needs to select no less than $1/\alpha$ scenarios, which provides that $N^{\min}=g(\epsilon, \alpha, n, r)=1/\alpha$ without the need of any PD function. Of course, this is just a simple intuitive explanation. In reality, $N^{\min}$ is also determined by other parameters, such as $\epsilon$ \cite{campi2008exact,calafiore2010random}. 
Observing from proposition 1, we conjecture that the analytical function $h$ is determined by the problem type (i.e. continuous or mixed-integer, and convex or nonconvex). Moreover, we think that $N^{\min}$ grows following $\epsilon$, $n$, and $r$ while it decreases as $\alpha$ increases.

\hfill$\square$

This hypothesis on the PD-free deterministication of (OU-ExU) is inspired by the existing findings on that of the (CCO-ExU) reviewed in the previous subsection. We provide a hypothesis rather than a proved theorem here since we found that the $N^{\min}$s determined by (\ref{N1})-(\ref{N3}) are unnecessarily large. It is worth pointing out that the conditions in proposition 1 are sufficient but not necessary for the $\epsilon$-optimality. In reality, engineers are generally interested in how to use as small numbers of scenarios as possible to achieve a high $\epsilon$ more than the exact expressions of $N^{\min}$. The relations in (\ref{hypo}) indicate that $N^{\min}$ grows following $\epsilon$, $n$, and $r$ while it decreases as $\alpha$ increases. Moreover, the $N^{\min}$s calculated in proposition 1 (i.e. the existing research) are unacceptably big for large-scale engineering systems, such as smart grids In next section, we will show that it is possible to achieve a high $\epsilon$ with a much smaller number of input scenarios $N$, i.e. $N \ll N^{\min}$. 

\section{Strategic Scenario Selection}
\label{sec:2}

\subsection{Theories of Active Scenarios and Repeated Scenarios}
In the theory of constrained optimization, an active constraint means the constraints that cause the limitation on the objective function \cite{chong2004introduction}. In other words, the optimal solution $x^*$ will not be changed by removing inactive constraints. Motivated by this characteristics of constrained optimization, we have the following definition and lemma.
\smallskip

\noindent
\textbf{Definition.} Let $\mathcal{N}$ be the set of $N$ scenarios \textit{randomly} selected from the historical measurements of uncertain variables $\mathcal{U}$, i.e. $\mathcal{N}=\{u^{(i)}\,(i=1,\ldots,N)\}$. For a given optimal solution $x^*$ of a GSP, $u^{(j)} \in \mathcal{N}$ is a active scenario if $g(x,u^{(j)})$ contain at lease one active constraints.
\smallskip

\noindent
\textbf{Lemma 1} (on active scenarios): \textit{For Case 1, the number of active scenarios $\Check{N}$ is less than or equal to the number $n$ of variables $x$ in} (CCO-ExU).
\smallskip

\textit{Proof}: According to Lemma 2.2 of \cite{calafiore2010random}, any finite dimensional (GSP-CCO-ExU) of Case 1 has at most $n$ active constraints if it is feasible. Since $g(\cdot)$ in (\ref{constr_CRP})/(\ref{constr_CSP}) is a vector of $m$ functions, each scenario $u^{(i)}$ contributes $m$ constraints to the optimization model. Even for the most conservative case, where each active scenario only contain one active constraint, the number of active scenarios $\Check{N}$ is at most $n$. As a result, $\Check{N}$ will not exceed $n$ for all other cases. 

\hfill$\square$ 

Lemma 1 and proposition 1 together imply $\Check{N} \ll N^{\min} \le N$ for Case 1 of (CCO-ExU). This relation is also applicable to nonconvex (OU-ExU) with integer decision variables, although the analytical expressions of $N^{\min}$ for these problems are still not exactly known. Moreover, we have the following proposition.
\smallskip

\noindent
\textbf{Proposition 2} (on repeated scenarios): \textit{There exist repeated scenarios in $\mathcal{N}$ which can be removed without impacting the solutions of GSPs.}
\smallskip

\textit{Proof}: In statistics, if the probability of a scenario $u^{(s)}$ is $\mathbb{P}[u^{(s)}]$, it will likely appear ($N*\mathbb{P}[u^{(s)}]$) times in set $\mathcal{N}$. However, in optimization, one $u^{(s)}$ is sufficient instead of ($N*\mathbb{P}[u^{(s)}]$), which means that the rest ($N*\mathbb{P}[u^{(s)}]$-1) $u^{(s)}$'s can be removed without impacting the optimal solution.

\hfill$\square$ 

\subsection{Methods of Strategic Scenario Selection}
 The purpose of strategic scenario selection (SSS) is to find $\tilde{N}$ ($\Check{N} \le \tilde{N} \ll N^{\min}$) scenarios that contains all active scenarios and $\tilde{N}$ is as small as possible. Inspired by proposition 2, we use the following dissimilarity logic to avoid repeated scenarios
 \begin{equation} \label{dissimilar}
     u^{(k+1)} \neq u^{(i)},\;(i=1,\ldots,k)
 \end{equation}
 which means the new selected scenario $u^{(k+1)} \in \mathcal{N}$ is different from all previous selected scenarios. It is straightforward to know that the logic (\ref{dissimilar}) can effectively eliminate repeated scenarios in $\mathcal{N}$. A numerical experiment showing the effectiveness of (\ref{dissimilar}) is given in Section 5.
 
\begin{figure}[h]
\centering
\includegraphics[width=0.8\textwidth]{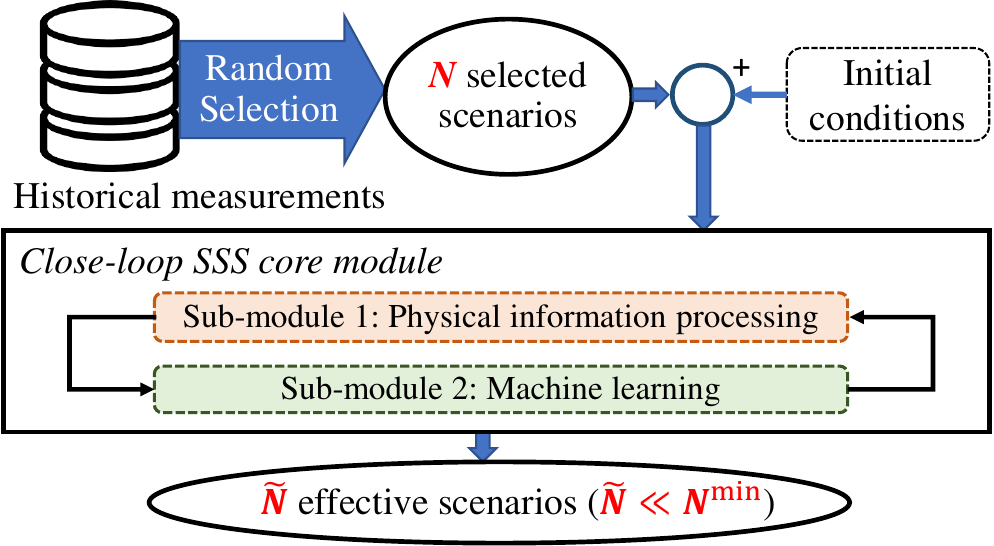}
 \caption{Flow chart of the close-loop learning-aided strategic scenario selection.}
  \label{fig:SSS}
\end{figure}

 In this research, we are interested in removing not only the repeated scenarios but also more inactive scenarios to further reduce the size of set $\mathcal{N}$. For an engineering perspective, we propose a close-loop learning-aided framework of SSS as shown in Figure \ref{fig:SSS}. First, many engineering systems like smart grids \cite{jenkins2012smart} are significantly impacted by some physical conditions, such as time span, weather, ambient temperature, and seasons. This physical information will be collected and processed in the first sub-module of the SSS core module. Then, in the second sub-module, properly selected machine learning algorithms which, on one hand, helps figure out the best physical information in sub-module 1 and, on the other hand, select effective scenarios guided by the physical information analyzed in sub-module 1. Here, we only provide introductory information about this advanced SSS framework while we'll provide a detailed discussion on it in a future engineering paper.
\smallskip

 \noindent
\textbf{Remark 1.} Although the hypothesis does not provide exact mathematical expressions on the needed number of scenarios for a desired confidence level, it offers a unique opportunity of effectively solving optimization under uncertainty for complex engineering systems. First, one can create a large scenario set $\mathcal{N}$. Then, under the help of the SSS methods, the numbers of scenarios that are finally input to (GSP-ExU) will be acceptably small with the same confidence level as inputting the original large set $\mathcal{N}$.

\section{Endogenous Uncertainty-compatible PD-free GSP}
This section presents another core contribution of this paper, i.e. the EnU- compatible PD-free GSP. This section will develop two paradigms for determinisfying (OU) into its deterministic equivalents in the form of GSP.

\subsection{A Non-iterative Paradigm of Deterministication}
We develop the following PD-free GSP as the deterministic equivalent of (OU):
\begin{subequations} \label{GSPEnU}
\begin{align}
\text{(GSP)}\quad\quad \min_{x} \quad & f_0(x) +\frac{1}{N}\sum_i^N (f_1(x,w^{(i)})+ \lambda_i (y_i-y^{(i)}))   \\
 \mathrm{s.t.} \quad  &g(x,w^{(i)}) \le 0,\;   (i=1,\ldots,N)  \\
  &y_i=q(v^{(i)},x),\quad ((w^{(i)},y^{(i)})\in \mathcal{W})  \label{JPD}
\end{align} 
\end{subequations}
where $\lambda_i$ is a Lagrangian multiplier, $y_i$ denotes the joint probability of scenario $w^{(i)}$ which is a function of $x$, and $q(v,x)=\prod_{j=1}^{r_2}q_j(v_j,x)$. The historical data set $\mathcal{U}$ is extended by including the joint probability $y^{(i)}$ of each scenario $w^{(i)}$ and denoted as $\mathcal{W}$. If function $q(\cdot)$ in (\ref{JPD}) is not given, both $y^{(i)}$ and $q$ can be obtained by training or fitting the historical data. Let $x^*_N$ be the optimal solution of (GSP) with $N \ge N^{\min}$ input scenarios, we have the following theorem.
\smallskip

\noindent
\textbf{Theorem 1} (on PD-free GSP under EnU): \textit{If the hypothesis holds true, then:}

\noindent
i. $x^*_N$ \textit{is a lower bound to the solution of} (OU) \textit{with a confidence level of} $\epsilon$.

\noindent
ii. $x^*_N$ \textit{is optimal to} (OU) \textit{with a confidence level of $\epsilon$ if} (OU) \textit{is convex and $q$ is affine on} $x$.
\smallskip

\textit{Proof}: First, we obtain the following re-formulation of (OU), which is in the form of optimization under ExU, by explicitly including the PD functions (\ref{PDn}) of EnU as a constraint: 
\begin{align}
 \min_{x} \quad & \text{(\ref{obj_OU})}  \nonumber \\
 \mathrm{s.t.} \quad  & \text{(\ref{constr_OU})} \label{reform-OU} \\
 & y = q(v,x)=\prod_{j=1}^{r_2}q_j(v_j,x)=\prod_{j=1}^{r_2}y_j(v_j) \nonumber
\end{align}
where the uncertain variables $u \in \mathbb{R}^{r_1}$ follow PD functions (\ref{PDx}) while the uncertain variables $v \in \mathbb{R}^{r_2}$ follow the PD functions $\mathbb{P}[v_j]=y_j(v_j)$ that is in the form of ExU PD functions.

Then, we consider a Lagrangian relaxation of (\ref{reform-OU}) :
\begin{align}
 \min_{x} \quad & \text{(\ref{obj_OU})} - \lambda(y-q(v,x)) \label{L-OnU} \\
 \mathrm{s.t.} \quad  & \text{(\ref{constr_OU})}, \nonumber
\end{align}
which is an (OU-ExU). According to (GSP-ExU), model (\ref{GSPEnU}) is exactly the deterministic GSP equivalent of (\ref{L-OnU}). Assuming the hypothesis true, the optimal solution $x^*_N$ of (\ref{GSPEnU}) is optimal to (\ref{L-OnU}) with a confidence level of $\epsilon$. It suffices to show that the optimal solution of (\ref{L-OnU}) is a lower bound of the solution of (\ref{reform-OU}), i.e. the solution of (OU), since (\ref{L-OnU}) is a Lagrangian relaxation of (\ref{reform-OU}). Further, the optimal solution of (\ref{L-OnU}) is also optimal to (\ref{reform-OU}) if (\ref{reform-OU}) is convex according to the zero-gap property of Lagrangian relaxation for convex problems. Model (\ref{reform-OU}) is convex if (OU) is convex with a $q$ function that is affine on $x$.

\hfill$\square$

Further, we consider a subset of (OU):
\begin{subequations} \label{CCO}
\begin{align}
\text{(CCO)}\quad\quad \min_{x} \quad & f_0(x) \label{obj_CCO} \\
 \mathrm{s.t.} \quad  & \mathbb{P}_{w}[g(x,w) \le 0] \ge 1-\alpha, \label{constr_CCO} 
\end{align}
\end{subequations}
which is a CCO under both ExU and EnU. Following (GSP), the corresponding deterministic equivalent of (CCO) in the form of GSP is given as:
\begin{subequations} \label{CCO-GSP}
\begin{align}
\text{(GSP-CCO)}\quad\quad \min_{x} \quad & f_0(x) +\frac{1}{N}\sum_i^N \lambda_i (y_i-y^{(i)})   \\
 \mathrm{s.t.} \quad  &g(x,w^{(i)}) \le 0,\;   (i=1,\ldots,N)  \\
  &y_i=q(v^{(i)},x),\quad ((w^{(i)},y^{(i)})\in \mathcal{W})
\end{align} 
\end{subequations}
 Let $x^*_N$ be the optimal solution of (GSP-CCO) with $N$ input scenarios, we have the following lemma of the above theorem.
 \smallskip
 
 \noindent
\textbf{Lemma 2} (on PD-free GSP for (CCO) under both ExU and EnU): 

\noindent
i. \textit{For Case 1,} $x^*_N$ \textit{is also optimal to} (CCO) \textit{with a confidence level of} $\epsilon$ \textit{if N satisfies condition} (\ref{N1}).

\noindent
ii. \textit{For Case 2,} $x^*_N$ \textit{is a lower bound to the solution of} (CCO) \textit{with a confidence level of} $\epsilon$ \textit{if N satisfies condition} (\ref{N2}).

\noindent
iii. \textit{For Case 3,} $x^*_N$ \textit{is a lower bound to the solution of} (CCO) \textit{with a confidence level of} $\epsilon$ \textit{if N satisfies condition} (\ref{N3}).
\smallskip
 
\textit{Proof}: According to proposition 1, the hypothesis holds true for these cases. It suffices to prove this lemma by applying the conclusions in the theorem. 
 
 \hfill$\square$
 
 \noindent
\textbf{Remark 2.} Lemma 2 indicates that, using (GSP) as the deterministic equivalent to solve an (OU) problem, the type of uncertain variables does not bring any difference to the needed number of scenarios.

\subsection{An Iterative Paradigm}

\begin{figure}[h]
\centering
\includegraphics[width=0.76\textwidth]{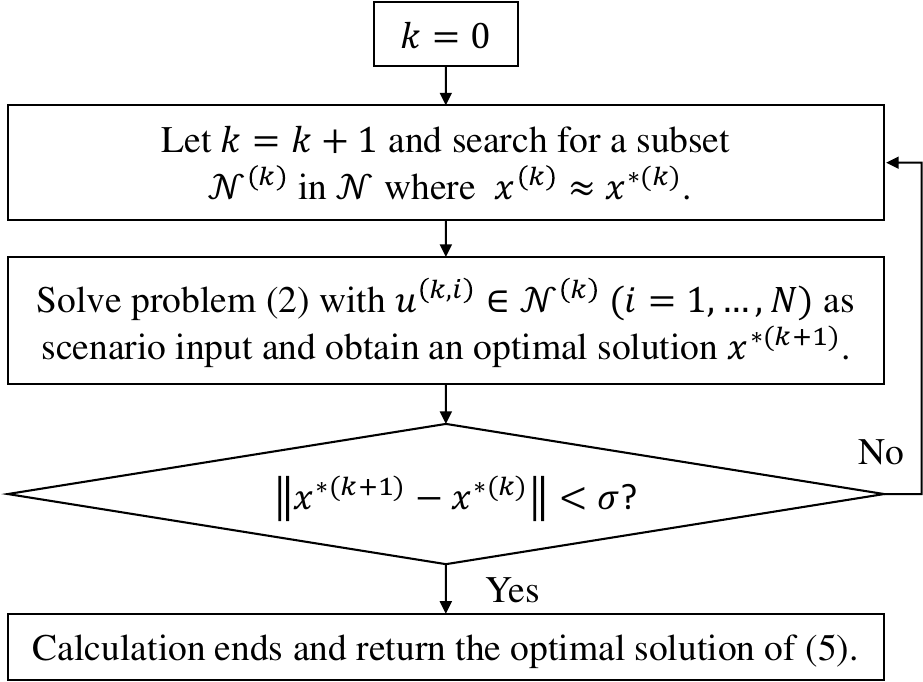}
  \caption{The iterative algorithm of endogenous uncertainty-compatible PD-free GSP.}
  \label{fig:flowchart}
\end{figure}

For determinisfying the optimization under EnU, the above non-iterative paradigm (GSP) preserves the beauty of (GSP-ExU). While many optimization problems in engineering systems are nonconvex, the optimality only holds for the convex cases under the above non-iterative framework. Moreover, the function $q$ in (\ref{JPD}) may be either unavailable or very expensive/time-consuming to obtain for some cases. As a remedy to the situation that the non-iterative paradigm can not work effectively, we propose an iterative paradigm as shown in Figure \ref{fig:flowchart}. For setting up this algorithm, $K*N$ scenarios are selected to produce a data set $\mathcal{N}=\{\mathcal{N}^{(1)},\ldots,\mathcal{N}^{(K)}\}$ where the $i$th data point in the $k$th subset is ($x^{(k)},w^{(k,i)}$). The advantage of this iterative paradigm resides in that it leverages the deterministic equivalent (GSP-ExU) for solving the optimization (OU).

 \section{Applications and Numerical Experiments}
This section presents two applications of the PD-free GSP solution paradigm in solving smart grid optimization (SGO) \cite{anjos2015optimization} problems. One is under exogenous uncertainty while the other is under endogenous uncertainty.

\subsection{Optimal Power Flow with High Renewable Energy}
This subsection considers a classical optimal power flow (OPF) \cite{momoh1997challenges} problem with high level of renewable energy as an example of SGO under exogenous uncertainty. Under the chance-constrained logical model, the OPF formulation is given as
\begin{subequations} \label{SG-CCO}
\begin{align}
 \min_{P^G,Q^G,\beta} \;\; &  \mathbb{E}_{u}[\sum_kf_{1,k}(P_k^G+\beta_k \sum_lu_l)] \label{obj_CCO} \\
 \mathrm{s.t.} \;\;  & \mathbb{P}_{u}[g(P^G, Q^G, \beta, u) \le \text{or} = 0] \ge 1-\alpha  \label{constr1_CCO} 
\end{align}
\end{subequations}
where the cost $c_k^{\prime}$ is a quadratic function of the real part of the real-time generation of bus $k \in$ bus set $\mathcal{B}$, decision variables $x=[P^G, Q^G,\beta]$, and uncertain variables $u$, that follow exogenous PD functions (\ref{PDx}), are the difference between the real-time (actual) renewable generation and the forecasted ones $P^R$. The detailed expressions of the vector function $g$ in (\ref{constr1_CCO}) are given as 
\begin{subequations} \label{ACPF}
\begin{align}
&  V^{\rm T}H^p_kV +P_k^L-P_k^R-u_k=(P_k^G+\beta_k\sum_lu_l) \label{active} \\
 &  V^{\rm T}H^q_kV +\eta(P_k^L-P_k^R- u_k)=(Q_k^G+\eta \beta_k\sum_lu_l) \label{reactive} \\
 & V^{\rm T}H_{ij}V \le \Bar{S}_{ij}^{Thermal} \label{thermal} \\
 & \underline{P}^G_k,\underline{Q}^G_k \le P^G_k, Q^G_k\le \Bar{P}^G_k, \Bar{Q}^G_k,\;\underline{V}_k\le \|V_k \|_2 \le \Bar{V}_k \label{bounds} 
\end{align}
\end{subequations}
where $P^G_k$ and $Q^G_k$ denote the bases of active and reactive power generation respectively, $P^L_k$ represents the active load, ($P_k^R-u_k$) is the real-time active power of renewable generation, and $\beta_k$ is the participation factor of the generator at bus $k$ \cite{jabr2013adjustable}.  For simplicity, we do not consider the load uncertainty in this application. Vector $V$ of the real and imaginary parts of bus voltages are state variables in this problem, while $V_k$ is a two-dimension sub-vector of $V$. Note that the state variables vary following the decision and control variables. The parameter $\eta$ implies that the loads and renewable generators are operated in the constant power factor mode. Constraint (\ref{thermal}) represents the thermal limits of power lines.

A deterministic equivalent of (\ref{SG-CCO}) based on the PD-free GSP model (\ref{GSP}) is established and given as 
\begin{subequations} \label{SG-GSP}
\begin{align}
 \min_{S^G,\beta} \;\;\; & \frac{1}{N}\sum_i^N \sum_kf_{1,k}(P_k^G+\beta_k \sum_lu_l^{(i)}) \label{obj_sggsp} \\
 \mathrm{s.t.} \;\;\;  &g(P^G, Q^G, \beta, u^{(i)}) \le \text{or} = 0,\; (i=1,\ldots,\tilde{N}) \label{constr1_sggsp} 
\end{align} 
\end{subequations}
where $u^{(i)} \in \tilde{\mathcal{N}}$ which is the set of scenarios selected by the SSS algorithm, and $\tilde{N}=|\tilde{\mathcal{N}}|$.

\subsection{OPF in Distribution Systems under Market-induced\\ Consumer Demand Uncertainty}
 Market-enabling is one of the key features of smart grids, where the locational marginal prices (LMP) of distribution nodes are determined by the optimal solution of distribution OPF \cite{sotkiewicz2006nodal,shaloudegi2012novel} while customers decide their electricity consumption based on the real-time LMPs. On one hand, the forecasted loads (i.e. consumer demands) are parameters for obtaining the optimal solution of decision variables. On the other hand, the energy prices in turn impacts the real-time loads while one can consider that the energy price is determined by the optimal solution of the DOPF. 

A distribution OPF problem under endogenous uncertainty (DOPF-EnU) is considered in this subsection, where the EnU comes from the consumer demands that depend on the energy prices. We use the following chance-constrained optimization formulation (\ref{SG-CCO1}) to capture the logical model of DOPF-EnU:
\begin{subequations} \label{SG-CCO1}
\begin{align}
 \min_{P^G,Q^G,\beta} \;\; &  \mathbb{E}_{v}[\sum_kf_{1,k}(P_k^G+\beta_k \sum_lv_l)] \label{obj_CCO1} \\
 \mathrm{s.t.} \;\;  & \mathbb{P}_{v}[g(P^G, Q^G, \beta, v) \le \text{or} = 0] \ge 1-\alpha  \label{constr1_CCO1} 
\end{align}
\end{subequations}
where uncertain variables $v$, that follow the endogenous PD functions (\ref{PDn}), are the difference between the real-time (actual) loads and the forecasted loads $P^L$. The detailed expressions of the vector function $g$ include (\ref{thermal}), (\ref{bounds}), and
\begin{subequations} \label{ACPF1}
\begin{align}
&  V^{\rm T}H^p_kV +(P_k^L+v_k)-(P_k^G+\beta_k\sum_lv_l)=0 \label{active} \\
 &  V^{\rm T}H^q_kV +\eta(P_k^L+ v_k)-(Q_k^G+\eta \beta_k\sum_lv_l)=0, \label{reactive} 
\end{align}
\end{subequations}
where ($P_k^L+v_k$) is the real-time active load at bus $k$. A deterministic equivalent of (\ref{SG-CCO1}) based on the PD-free GSP model (\ref{GSPEnU}) is created and given as follows
\begin{subequations} \label{SG-GSP1}
\begin{align}
 \min_{S^G,\beta} \;\;\; & \frac{1}{N}\sum_i^N \{ \sum_kf_{1,k}(P^G_k+\beta_k \sum_lv_l^{(i)}) + \lambda_i(y_i-y^{(i)})\} \label{obj_sggsp1} \\
 \mathrm{s.t.} \;\;\;  & g(P^G, Q^G, \beta, v^{(i)}) \le \text{or} = 0, \;(i=1,\ldots,\tilde{N}) \nonumber \\
 & y_i=q(v^{(i)},P^G, Q^G,\beta) \label{constr1_sggsp1}
\end{align} 
\end{subequations}
where $(v,y)^{(i)} \in \tilde{\mathcal{V}}$.

\subsection{Numerical Experiments}
\subsubsection{PD-free v.s. PD-based}

\begin{table}[h]
\centering
\caption{Needed scenarios of the PD-based and PD-free deterministication methods for the power test cases.}
\begin{tabular}{|c|c|c|c|c|}
\hline
\multirow{2}{*}{\textbf{\begin{tabular}[c]{@{}c@{}}IEEE system\end{tabular}}} & \multirow{2}{*}{\textit{\textbf{n}}} & \multirow{2}{*}{\textit{\textbf{r}}} & \multicolumn{2}{c|}{\textbf{Needed of scenarios}} \\ \cline{4-5} 
                                                                                &                                      &                                      & \textbf{PD-based}        & \textbf{PD-free}      \\ \hline
\textbf{9-bus}                                                                  & 22                                    & 3                                    & 125                        & 768                      \\ \hline
\textbf{57-bus}                                                                 & 124                                    & 7                                   & 78k                      & 4k                     \\ \hline
\textbf{118-bus}                                                                & 305                                   & 19                                   & 1.9*$10^{13}$                      & 9.8k                    \\ \hline
\end{tabular}
\end{table}

This experiment evaluates the performance of the PD-free deterministication via comparing it with that of PD-based deterministication on problem (\ref{SG-CCO}). We uses the convex relations of constraints (\ref{ACPF}) \cite{li2016convex,li2016non,li2017convex,li2016theconvex} in this experiment. As a result, the problem (\ref{SG-CCO}) is convex so that we can use expression (\ref{N1}) to calculate the needed numbers of scenarios determined by the PD-free methods. The needed numbers of the PD-based and PD-free scenarios on some IEEE test systems are tabulated in the Table 2. In this numerical experiment, we considered a constraint violation probability of $\alpha=5\%$ and an accuracy of 90$\%$ (i.e. $\epsilon=90\%$). For the PD-based sampling, we assume that all ExU variables follow the normal distributions and 5 samples per normal distribution can provide 90$\%$ accuracy. It can be observed from the table that the needed number of PD-based scenarios increases exponentially over the number of uncertain variable. The needed number of PD-free scenarios is much less than that of the PD-based ones for large systems.

\subsubsection{With SSS v.s. Without SSS}
\begin{figure}
    \centering
    \includegraphics[width=0.6\textwidth]{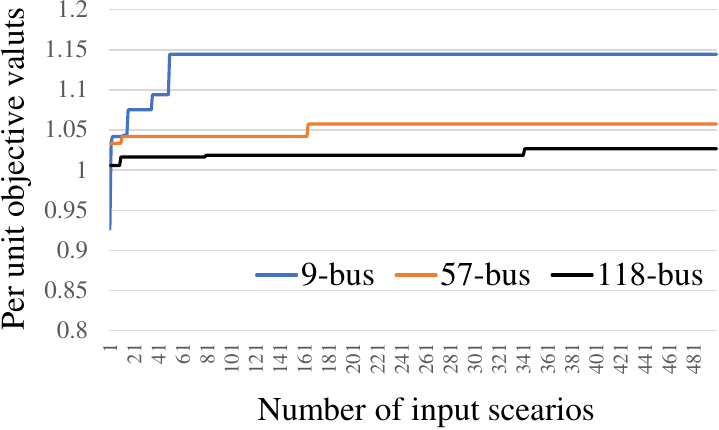}
  \caption{Numerical results of dissimilarity -based scenario selection. Each step-up on the curves means an active scenario has been added.}
  \label{fig:SSSresult}
\end{figure}

This numerical experiment solves an optimal multi-period power system operation problem with energy storage under the uncertainty of both loads and renewable energy resources whose one single period formulation is highly similar to (\ref{SG-CCO}). In other words, the model size of this problem is multiple times as (\ref{SG-CCO})'s. The effectiveness of the dissimilarity logic (7) is evaluated through solving such a computationally challenging problem on the IEEE 9-, 57- and 118-bus systems, where an operation horizon of 24 hours is considered. In this experiment, scenarios are added to the PD-free GSP (\ref{GSP}) one-by-one of which the simulation results are plotted in Figure \ref{fig:SSSresult}. For the 9-bus case, after adding around 50 scenarios selected by the dissimilarity-based method, the objective value no longer increases no matter how many scenarios will be added. That means the solution process has reached the optimal solution. These numbers for the 57- and 118-bus cases are around 160 and 340 respectively. It is worth noting that, if we consider the constraints in (\ref{ACPF}), $N^{\min}_{9}$=44k, $N^{\min}_{57}$=236k and $N^{\min}_{118}$=693k according to equation (\ref{N1}). In other words, the  dissimilarity-based scenario selection dramatically reduces the needed number of scenarios by effectively removing the repeating scenarios. In our future research, we will seek to develop more advanced SSS methods for more effective scenario reduction as described in Subsection 3.2.

\section{Conclusion}

This paper establishes a novel solution paradigm--the probability distribution (PD)-free general scenario programming--for general optimization under both exogenous and endogenous uncertainties based on an existing PD-free deterministication methods. According to the discussions throughout the paper, we can conclude that our contributions includes: 1) generalizing the PD-free deterministication method, which originally works for chance-constrained optimization only, and makes it applicable to more general optimization problems under uncertainty; 2) extending the PD-free deterministication method, which originally works for problems under exogenous uncertainty only, and makes it applicable to problems under endogenous uncertainty; 3) revealing more properties of the PD-free deterministication, such as the properties on active scenarios and repeated scenarios; and 4) developing strategic scenarios selection methods which can effective reduce the needed scenarios for a desired confidence level. The applications of the proposed approaches on two smart grid optimization problems under exogenous and endogenous uncertainties shows that they have a high capability of solving related optimization problems of complex engineering systems like smart grids.

%
%

\bibliographystyle{unsrt}
\bibliography{main}

\begin{thebibliography}{10}

\bibitem{nohadani2018optimization}
Omid Nohadani and Kartikey Sharma.
\newblock Optimization under decision-dependent uncertainty.
\newblock {\em SIAM Journal on Optimization}, 28(2):1773--1795, 2018.

\bibitem{li2008chance}
Pu~Li, Harvey Arellano-Garcia, and G{\"u}nter Wozny.
\newblock Chance constrained programming approach to process optimization under
  uncertainty.
\newblock {\em Computers \& chemical engineering}, 32(1-2):25--45, 2008.

\bibitem{birge2011introduction}
John~R Birge and Francois Louveaux.
\newblock {\em Introduction to stochastic programming}.
\newblock Springer Science \& Business Media, 2011.

\bibitem{bertsimas2011theory}
Dimitris Bertsimas, David~B Brown, and Constantine Caramanis.
\newblock Theory and applications of robust optimization.
\newblock {\em SIAM review}, 53(3):464--501, 2011.

\bibitem{nocedal2006numerical}
Jorge Nocedal and Stephen Wright.
\newblock {\em Numerical optimization}.
\newblock Springer Science \& Business Media, 2006.

\bibitem{anand2017comparative}
Rimmi Anand, Divya Aggarwal, and Vijay Kumar.
\newblock A comparative analysis of optimization solvers.
\newblock {\em Journal of Statistics and Management Systems}, 20(4):623--635,
  2017.

\bibitem{campi2008exact}
Marco~C Campi and Simone Garatti.
\newblock The exact feasibility of randomized solutions of uncertain convex
  programs.
\newblock {\em SIAM Journal on Optimization}, 19(3):1211--1230, 2008.

\bibitem{calafiore2010random}
Giuseppe~Carlo Calafiore.
\newblock Random convex programs.
\newblock {\em SIAM Journal on Optimization}, 20(6):3427--3464, 2010.

\bibitem{campi2018general}
Marco~Claudio Campi, Simone Garatti, and Federico~Alessandro Ramponi.
\newblock A general scenario theory for nonconvex optimization and decision
  making.
\newblock {\em IEEE Transactions on Automatic Control}, 63(12):4067--4078,
  2018.

\bibitem{calafiore2012mixed}
Giuseppe~C Calafiore, Daniel Lyons, and Lorenzo Fagiano.
\newblock On mixed-integer random convex programs.
\newblock In {\em 2012 IEEE 51st IEEE Conference on Decision and Control
  (CDC)}, pages 3508--3513. IEEE, 2012.

\bibitem{esfahani2014performance}
Peyman~Mohajerin Esfahani, Tobias Sutter, and John Lygeros.
\newblock Performance bounds for the scenario approach and an extension to a
  class of non-convex programs.
\newblock {\em IEEE Transactions on Automatic Control}, 60(1):46--58, 2014.

\bibitem{chong2004introduction}
Edwin~KP Chong and Stanislaw~H Zak.
\newblock {\em An introduction to optimization}.
\newblock John Wiley \& Sons, 2004.

\bibitem{jenkins2012smart}
Nick Jenkins, Kithsiri Liyanage, Jianzhong Wu, and Akihiko Yokoyama.
\newblock {\em Smart Grid}.
\newblock Wiley., 2012.

\bibitem{anjos2015optimization}
Miguel~F Anjos.
\newblock Optimization for power systems and the smart grid.
\newblock In {\em Modeling and Optimization: Theory and Applications}, pages
  29--47. Springer, 2015.

\bibitem{momoh1997challenges}
JA~Momoh, RJ~Koessler, MS~Bond, B~Stott, D~Sun, A~Papalexopoulos, and
  P~Ristanovic.
\newblock Challenges to optimal power flow.
\newblock {\em IEEE Transactions on Power systems}, 12(1):444--455, 1997.

\bibitem{jabr2013adjustable}
Rabih~A Jabr.
\newblock Adjustable robust opf with renewable energy sources.
\newblock {\em IEEE Transactions on Power Systems}, 28(4):4742--4751, 2013.

\bibitem{sotkiewicz2006nodal}
Paul~M Sotkiewicz and Jesus~M Vignolo.
\newblock Nodal pricing for distribution networks: efficient pricing for
  efficiency enhancing dg.
\newblock {\em IEEE transactions on power systems}, 21(2):1013--1014, 2006.

\bibitem{shaloudegi2012novel}
Kiarash Shaloudegi, Nazli Madinehi, SH~Hosseinian, and Hossein~Askarian
  Abyaneh.
\newblock A novel policy for locational marginal price calculation in
  distribution systems based on loss reduction allocation using game theory.
\newblock {\em IEEE transactions on power systems}, 27(2):811--820, 2012.

\bibitem{li2016convex}
Qifeng Li, Raja Ayyanar, and Vijay Vittal.
\newblock Convex optimization for des planning and operation in radial
  distribution systems with high penetration of photovoltaic resources.
\newblock {\em IEEE Transactions on Sustainable Energy}, 7(3):985--995, 2016.

\bibitem{li2016non}
Qifeng Li and Vijay Vittal.
\newblock Non-iterative enhanced sdp relaxations for optimal scheduling of
  distributed energy storage in distribution systems.
\newblock {\em IEEE Transactions on Power Systems}, 32(3):1721--1732, 2016.

\bibitem{li2017convex}
Qifeng Li and Vijay Vittal.
\newblock Convex hull of the quadratic branch ac power flow equations and its
  application in radial distribution networks.
\newblock {\em IEEE Transactions on Power Systems}, 33(1):839--850, 2017.

\bibitem{li2016theconvex}
Qifeng Li and Vijay Vittal.
\newblock The convex hull of the ac power flow equations in rectangular
  coordinates.
\newblock In {\em 2016 IEEE Power and Energy Society General Meeting (PESGM)},
  pages 1--5. IEEE, 2016.

\end{thebibliography}

\end{document}